\documentclass[12pt]{amsart}
\usepackage{amsmath}
\usepackage{amsfonts,amssymb,amsthm}
\usepackage{graphicx}

\newcommand{\beq}{\begin{equation}}
\newcommand{\eeq}{\end{equation}}
\newcommand{\bbar}{\begin{eqnarray}}
\newcommand{\eear}{\end{eqnarray}}

\newcommand{\thm}[2]{\begin{#1} #2 \end{#1}}

\newtheorem{theorem}{Theorem}[section]
\newtheorem{fact}[theorem]{Fact}
\newtheorem{itheorem}{Theorem}[section]
\newtheorem{lemma}[theorem]{Lemma}
\newtheorem{ilemma}[itheorem]{Lemma}

\begin{document}

\title[Harmonic Mean, \textit{etc.}]{Harmonic mean, random polynomials and stochastic matrices}

\author{Natalia L. Komarova}
\address{Institute For Advanced Study, Einstein Drive, Princeton, NJ
$08540$}
\address{Department of Applied Mathematics, University of
Leeds, Leeds LS2 9JT, UK}

\email{natalia@ias.edu}
\author{Igor Rivin}
\address{Mathematics department, University of Manchester,
Oxford Road, Manchester M13 9PL, UK}
\address{Mathematics Department, Temple University,
Philadelphia,  PA 19122 }
\address{Mathematics Department, Princeton University, Princeton,
NJ 08544}
\email{irivin@math.princeton.edu} \thanks{N.Komarova gratefully
acknowledges support from the Packard Foundation, the Leon Levy and
Shelby White Initiatives Fund, the Florence Gould Foundation, the
Ambrose Monell Foundation, the Alfred P. Sloan Foundation and the NSF.
I. Rivin would like to thank the EPSRC and the NSF for support}

\begin{abstract}
Motivated by a problem in learning theory, we are led to study the
dominant eigenvalue of a class of random matrices. This turns out to
be related to the roots of the derivative of random polynomials
(generated by picking their roots uniformly at random in the interval
$[0, 1]$, although our results extend to other distributions). This,
in turn, requires the study of the statistical behavior of the
harmonic mean of random variables as above, and that, in turn, leads
us to delicate question of the rate of convergence to stable laws and
tail estimates for stable laws.
\end{abstract}
\maketitle

\renewcommand{\theitheorem}{\Alph{itheorem}}
\section*{Introduction}
The original motivation for the work in this paper was provided by the
first-named author's research in learning theory, specifically in
various models of language acquisition (see \cite{knn,nkn,kn}) and
more specifically yet by the analysis of the speed of convergence of
the \emph{memoryless learner algorithm}. The setup is described in
some detail in Section \ref{memoryless}; here we will just recall
the essentials. There is a collection of concepts $R_1, \dots, R_n$
and words which refer to these concepts, sometimes ambiguously. The
teacher generates a stream of words, referring to the concept
$R_1$. This is not known to the student, but he must learn by, at each
steps, guessing some concept $R_i$ and checking for consistency with
the teacher's input.  The memoryless learner algorithm consists of
picking a concept $R_i$ at random, and sticking by this choice, until
it is proven wrong. At this point another concept is picked randomly,
and the procedure repeats. It is clear that once the student hits on
the right answer $R_1$, this will be his final answer, so the question
is then:
\begin{quote}
\emph{How quickly does this method converge to the truth?}
\end{quote}
Since the method is memoryless, as the name implies, it is clear that
the learning process is a \emph{Markov chain}, and as is well-known
the convergence rate is determined by the gap between the top
(Perron-Frobenius) eigenvalue and the second largest
eigenvalue. However, we are also interested in a kind of a
\emph{generic} behavior, so we assume that the sizes of overlaps
between concepts are \emph{random}, with some (sufficiently regular)
probability density function supported in the interval $[0, 1]$, and
that the number of concepts is large. This makes the transition matrix
random (that is, the entries are random variables) -- the precise
model is described in Section \ref{memoryless}. The analysis of
convergence speed then comes down to a detailed analysis of the size
of the second largest eigenvalue and also of the properties of the
eigenspace decomposition (the contents of Section \ref{eig}.) Our main
results for the original problem (which is presented in Section
\ref{convmem}) can be summarized in the following
\footnote{here and throughout this section, the reference to the relevant
theorem (lemma) in the main body of the paper is given in square
brackets}:
\begin{itheorem}{\bf [\ref{main}, \ref{main1}]}
Let $N_\Delta$ be the number of steps it takes for the student 
to have probability $1 - \Delta$ of learning the
concept. Then we have the following estimates for $N_\Delta$:
\begin{itemize}
\item
if the distribution of overlaps is \emph{uniform}, or more 
generally, the density function $f(1-x)$  at $0$ has the form 
$f(x) = c + O(x^\delta),$ $\delta, c > 0,$ then there exist positive
constants $C_1, C_2$ such that 
$$\lim_{n \rightarrow \infty} \mathbf{P}\left(C_1 <
\frac{N_\Delta}{|\log \Delta|n \log n} < C_2\right) = 1,$$

\item 
if the probability density function $f(1-x)$ is asymptotic to $c x^\beta
+ O(x^{\beta + \delta}), \quad \delta, \beta > 0$, as $x$ approaches
$0$, then we have 
$$\lim_{n \rightarrow \infty} \mathbf{P}\left(C'_1 <
\frac{N_\Delta}{|\log \Delta|n} < C'_2\right) = 1$$
for some positive constants $C'_1$ and $C_2'$, 
\item 
if the asymptotic behavior is as above, but $-1 < \beta < 0$, then
$$\lim_{x \rightarrow \infty}  \mathbf{P}\left(\frac{1}{x} < \frac{N_\Delta}{|\log
\Delta| n^{1/(1+\beta)}} < x\right) = 1.$$
\end{itemize}
\end{itheorem}

It should be said that our methods give
quite precise estimates on the constants in the asymptotic estimate,
but the rate of convergence is rather poor -- logarithmic -- so these
precise bounds are of limited practical importance.

\subsection*{Notation}
We shall use the notation $a \asymp b$ to mean that $a$ is
asymptotically the same as $b$. We say that $a \sim b$ if $a$ and $b$
have the same order of growth (in other words, there exist constants
$c_1$, $c_2$, $d_1$, $d_2$, with $c_1,c_2>0$, so that $c_1 a + d_1
\leq b \leq c_2 a + d_2.$) In addition we denote the expectation of a
random variable $x$ by $\mathbf{E}(x)$.

\subsection*{Eigenvalues and polynomials} In order to calculate the
convergence rate of the learning algorithm described above, we need to
study the spectrum of a class of random matrices. The matrix $T =
(T_{ij})$ is an $n \times n$ matrix with entries
 (see Section \ref{memoryless}):
\begin{equation*}
T_{ij} = 
  \begin{cases}
    a_i & i=j, \\
    \frac{(1-a_i)}{n-1} & \text{otherwise}.
  \end{cases}
\end{equation*}
Let $B=\frac{n-1}{n}(I-T)$, so that the eigenvalues of $T$,
$\lambda_i$, are related to the eigenvalues of $B$, $\mu_i$, by
$\lambda_i=1-n/(n-1)\mu_i$. In Section \ref{charderiv} we show the
following result:
\begin{ilemma}
\label{charpd}{\bf [\ref{amuz}]}  Let $p(x) = 
(x - x_1) \dots (x-x_n)$, where $x_i=1-a_i$. Then the characteristic
polynomial $p_B$ of $B$ satisfies:
$$
p_B(x) = \frac{x}{n} \frac{dp(x)}{dx}.
$$
\end{ilemma}

Lemma \ref{charpd}  brings us to the
following question:
\begin{quote}
\textbf{Question $1$:} Given a random polynomial $p(x)$ whose
roots are all real, and distributed in a prescribed way, what can
we say about the distribution of the roots of the derivative
$p^\prime(x)$?
\end{quote}

And more specifically, since the convergence behavior of $T^N$ is
controlled by the top of the spectrum: 

\begin{quote}
\textbf{Question $1^\prime$:} What can we say about the distribution
of the \emph{smallest} root of $p^\prime(x)$, given that the smallest
root of $p(x)$ is fixed?
\end{quote}

For Question $1^\prime$ we shall clamp the smallest root of $p(x)$ at
$0$. Letting $H_{n-1}$ be the \emph{harmonic mean} of the other
roots of $p(x)$ (which are all greater than zero with probability
$1$), our first observation will be 
\begin{ilemma}{\bf [\ref{bounds}]}
\label{ismall}
 The smallest root $\mu_*$ of $p^\prime(x)$
satisfies:
$$
\frac{1}{2}H_{n-1} \leq (n-1) \mu_* \leq H_{n-1}.
$$
\end{ilemma}
We will henceforth assume that the roots of the polynomial $p(x)$ are
a sample of size $\deg p(x)$ of a random variable, $x$, distributed in
the interval $[0, 1]$.  In this stochastic setting, it will be shown
that $(n-1) \mu_*$ tends to the harmonic mean of the non-zero roots of
$p$ with probability $1$, when $n$ is large. It then follows that the
study of the distribution of $\mu_*$ entails the study of the
asymptotic behavior of the harmonic mean of a sample drawn from a
distribution on $[0, 1]$.

\subsection*{Statistics of the harmonic mean.}
In view of the long and honorable history of the harmonic mean, it
seems surprising that its limiting behavior has not been studied more
extensively than it has. Such, however, does appear to be the case. It
should also be noted that the arithmetic, harmonic, and geometric
means are examples of the ``conjugate means'', given by
\begin{equation*}
m_{\mathcal F}(x_1, \dots, x_n) = {\mathcal F}^{-1} \left(\frac{1}{n} \sum_{i=1}^n 
{\mathcal F}(x_i)\right),
\end{equation*}
where ${\mathcal F}(x) = x$ for the arithmetic mean, ${\mathcal F}(x)
= \log(x)$ for the geometric mean, and ${\mathcal F}(x) = 1/x$ for the
harmonic mean. The interesting situation is when ${\mathcal F}$ has a
singularity in the support of the distribution of $x$, and this case
seems to have been studied very little, if at all. Here we will devote
ourselves to the study of harmonic mean. 

Given $x_1, \dots, x_n$ -- a sequence of independent, identically
distributed in $[0, 1]$ random variables (with common probability
density function $f$), the nonlinear nature of the harmonic mean leads
us to consider the random variable
\begin{equation*}
X_n = \frac{1}{n} \sum_{i=1}^n \frac{1}{x_i}. 
\end{equation*}
Since the variables $1/x_i$ are easily seen to have infinite
expectation and  variance, our prospects seem grim at first 
blush, but then we notice that the variable $1/x_i$ falls 
straight into the framework of the ``stable laws'' of L\'evy -- 
Khintchine (\cite{feller}). Stable laws are defined and discussed in
Section \ref{stablaw}. Which particular stable law comes up 
depends on the distribution function $f(x)$. If we assume that
\begin{equation*}
f(x) \asymp c x^\beta,
\end{equation*}
as $x \rightarrow 0$ (for the uniform distribution $\beta = 0, \quad c
= 1$), we have 
\begin{itheorem}
\label{stableob} If $\beta = 0$, then let $Y_n = X_n - \log n.$ 
The variables $Y_n$ converge in distribution to the unbalanced stable
law $G$ with exponent $\alpha=1$. If $\beta > 0$, then $X_n$ converges
in distribution to $\delta(x-{\mathcal E})$, where ${\mathcal E} =
\mathbf{E}(1/x),$ and $\delta$ denotes the Dirac delta function. If
$-1 < \beta < 0$, then $n^{1-1/(1+\beta)} X_n$ converges in
distribution to a stable law with exponent $\alpha=1+\beta$.
\end{itheorem}
The above result points us in the right direction, since it allows us
to guess the form of the following results ($H_n$ is the harmonic mean
of the variables):
\begin{itheorem}{\bf [\ref{log}, \ref{log1}]}
\label{harmlog}
 Let $H_n = 1/X_n$ and $\beta = 0$. Then there exists a constant
$\mathfrak{C}_1$ such that
\begin{equation*}
\lim_{n\rightarrow \infty} 
 \mathbf{E}(H_n\log n ) = \mathfrak{C}_1.
\end{equation*}
\end{itheorem}
\begin{itheorem}{\bf [\ref{exp}]}
\label{harmmu}
Suppose $\beta > 0$, let $y = 1/x$, and let ${\mathcal E}$ be the mean of the
variable $y.$ Then  
$$\mathbf\lim_{n\rightarrow
\infty} \mathbf{E}({\mathcal E} H_n) = 1.$$
\end{itheorem}
Finally,
\begin{itheorem}{\bf [\ref{negbeta}]}
\label{expneg}
Suppose $\beta < 0.$ Then there exists a constant $\mathfrak{C}_2$ such
that 
$$\mathbf{E}(H_n/n^{1-1/{(1+\beta)}}) = \mathfrak{C}_2.$$
\end{itheorem}
\begin{itheorem}[\ref{weaklaw},\ref{weakal}, {\rm Law of large numbers for 
harmonic mean}] \label{iweak} 
Let $\beta = 0$ and let $a  
> 0$. Then 
\begin{equation*}
\lim_{n\rightarrow \infty} \mathbf{P}( \vert H_n \log n -\mathfrak{C}_1\vert 
> a) = 0,
\end{equation*}
where $\mathfrak{C}_1$ is as in the statement of Theorem
\ref{harmlog}. If $\beta > 0$, and ${\mathcal E}$ is as in the
statement of Theorem \ref{harmmu}, then

\begin{equation*}
\lim_{n\rightarrow \infty} \mathbf{P}(\vert H_n 
-\frac{1}{{\mathcal E}}\vert) > a) = 0.
\end{equation*}
\end{itheorem}

The proofs of the results for $\beta = 0$ require estimates of the speed of
convergence in Theorem \ref{stableob}. The speed of convergence
results we obtain (in Section \ref{sconv}) are not best possible, but
the arguments are simple and general.  The estimates can be summarized
as follows:

\begin{itheorem}{\bf [\ref{convthm}]}
\label{convspeed}  Assume $\beta=0$. Let $g_n$ be the density
associated to $X_n- \log n$, and let $g$ be the probability density of
the unbalanced stable law with exponent $\alpha=1$.  Then we have
(uniformly in $x$):
\begin{equation*}
g_n(x) = g(x) + O(\log^2n/n).
\end{equation*}
\end{itheorem}

In addition to the laws of large numbers we have the following
limiting distribution results:
\begin{itheorem}{\bf [\ref{limlaw}, \ref{limlawxx}]}
\label{hlimlaw} For $\alpha=1$, the random variable $\log n(H_n\log n 
- {\mathfrak C}_1)$ converges to a variable with the distribution
function $1-G(-x/{\mathfrak C}_1^2),$ where $G$
is the limiting distribution (of exponent $\alpha=1$) of variables
$Y_n = X_n - c\log n$ and ${\mathfrak C}_1=1/c$.
\end{itheorem}

\begin{itheorem}{\bf [\ref{limlaw2}]}
\label{hlimlaw2}
For $\alpha>1$, the random variable $n^{1-1/\alpha} (H_n -
\frac{1}{{\mathcal E}})$ converges in distribution to a variable with
distribution function $1-G(-x {\mathcal
E}^2)$, where $G$ is the unbalanced stable distribution of  exponent $\alpha$.
\end{itheorem}

\begin{itheorem}{\bf [\ref{limlaw3}]}
\label{hlimlaw3}
For $0<\alpha<1$, the random variable $H_n/n^{1-1/\alpha}$ converges in
distribution to the variable with distribution function $1-G(1/x)$,
where $G$ is the unbalanced  stable distribution of exponent $\alpha.$
\end{itheorem}

The paper is organized as follows. In Section \ref{harme} we study
some statistical properties of a harmonic mean of $n$ variables and in
particular, find the expected value of its mean as $n\to \infty$. In
Section \ref{classpol} we explore the connection between the harmonic
mean and the smallest root of the derivative of certain random
polynomials. In Section \ref{classmat} we uncover the connection
between the rate of convergence of the memoryless learner algorithm,
eigenvalues of certain stochastic matrices and the harmonic mean. The
more technical material can be found in the Appendix. In Section
\ref{explic} of the Appendix we present an explicit derivation of the
stable law for a particular example with $\alpha=1$. 
In Section \ref{sconv} we evaluate the rate of convergence
of the distribution of the inverse of the harmonic mean to its stable
law.

\section{\label{harme}Harmonic mean}
\subsection{\label{prelim}Preliminaries}
Let $x_1, \dots, x_n$ be positive real numbers. The
\textit{harmonic mean}, $H_n$, is defined by
\begin{equation}
\frac{1}{H_n} = \frac{1}{n}\left(\sum_{i=1}^n \frac{1}{x_i}\right).
\end{equation}
Let $x_1, \dots, x_n$ be independent random variables, identically
uniformly distributed in $[0, 1]$. We will study statistical
properties of their harmonic mean, $H_n$, with emphasis on limiting
behavior as $n$ becomes large.

We will use auxiliary variables $X_n$ and $Y_n$, defined as
\begin{equation}
\label{XY}
X_n = \frac{1}{n}\sum_{i=1}^n \frac{1}{x_i} = \frac{1}{H_n},\,\,\,\,\,\,
Y_n=X_n-\log n,
\end{equation}
and also variables $y_i = \frac{1}{x_i}$. The distribution of $y_i$ is
easily seen to be given by
\begin{equation}
\label{defF}
\mathfrak{F}(z) = \mathbf{P}(y_i < z) = \begin{cases} 0 & z < 1\\ 1 -
\frac{1}{z} & \text{otherwise}. \end{cases}
\end{equation}

A quick check reveals that $y_i$ has infinite mean and variance, so
the Central Limit Theorem is not much help in the study of $X_n$.
Luckily, however, $X_n$ converges to a \emph{stable law}, as we 
shall see. A very brief introduction to stable laws is given in 
the next section.

\subsection{\label{stablaw}Stable limit laws}
Consider an infinite sequence of independent identically distributed
random variables $y_1, \dots, y_n, \dots$, with some probability
distribution function, $\mathfrak{F}$. Typical questions studied in
probability theory are the following.
\begin{quote}
Let $S_n = \sum_{j=1}^n y_j$. How is $S_n$ distributed? What can 
we say about the distribution of $S_n$ as $n \rightarrow \infty$?
\end{quote}
The best known example is one covered by the Central Limit Theorem of
de~Moivre - Laplace: if $\mathfrak{F}$ has finite mean ${\mathcal E}$ and
variance $\sigma^2$, then $(S_n - n {\mathcal E})/(\sqrt{n} \sigma)$ converges
in distribution to the \emph{normal distribution} (see, e.g.,
\cite{feller}).  In view of this result, one says that the variable
$X$ \emph{belongs to the domain of attraction of a non-singular
distribution $G$}, if there are constants $a_1, \dots, a_n, \dots$ and
$b_1, \dots, b_n, \dots$ such that the sequence of variables
$Y_n\equiv a_n S_n - b_n$ converges in distribution to $G$.  It was
shown by L\'evy and by Khintchine that having a domain of attraction
constitutes severe restrictions on the distribution as well as the
{\it norming sequences} $\{a_n\}$ and $\{b_n\}$. To wit, one can
always pick $a_n = n^{- 1/\alpha}, \quad 0 < \alpha \leq 2$.  
It turns out that $\alpha$ is determined by
the limiting behavior of the distribution $\mathfrak{F}$ so that
\begin{equation}
\label{oh1}
\lim_{|x|\to\infty}\frac{1}{|x|^{\alpha+1}}\frac{d \mathfrak{F}(x)}{dx}=\begin{cases} Cp, & x>0 \\ Cq, & x<0,\end{cases} 
\end{equation}
where $p+q=1$. In that case, $G$ is called a \emph{stable distribution
of exponent $\alpha$}. Note that the case when $\alpha = 2$
corresponds to the Central Limit Theorem. If the variable $y$ belongs
to the domain of attraction of a stable distribution of exponent $\alpha > 1$, then
$y$ has a finite expectation ${\mathcal E}$; just as in the case
$\alpha = 2$, we can choose $b_n = n^{1-1/\alpha} {\mathcal E}.$ When
$\alpha < 1$, the variable $y$ does not have a finite expectation, and
it turns out that we can take $b_n \equiv 0$; for $\alpha = 1$, we can
take $b_n = c \log n$, where $c$ is a constant depending on
$\mathfrak{F}$. Thus, the normal distribution is a stable distribution
of exponent $2$ (and it is also unique, up to scale and shift). This
is one of the few cases where we have an explicit expression for the
density of a stable distribution; in other cases we only have
expressions for their \emph{characteristic functions}. The
characteristic function $\Psi(k)$ of a distribution function $G(x)$ is
defined to be $\int_{-\infty}^\infty \exp(i k x) d G(x),$ that is, as
the Fourier transform of the density function. Levy and Khinchine
showed that the characteristic functions of stable distributions can
be parameterized as follows:
\begin{equation}
\label{expforms}
\begin{split}
\log  \Psi(k) &= 
\begin{cases} 
C \frac{\Gamma(3-\alpha)}{\alpha (\alpha - 1)} 
\left[\cos \frac{\pi \alpha}{2} -{\rm sign}\, k i (p-q) \sin\frac{\pi \alpha}{2} 
\right]|k|^\alpha, & \alpha \neq 1\\
- C\left[\frac{1}{2} \pi -{\rm sign}\, k i (p-q) \log |k|\right] 
|k| + const,  & \mbox{$\alpha = 1$},
\end{cases}
\end{split}
\end{equation}
where the constants $C$, $p$ and $q$ can be defined by the following
limits:
\begin{eqnarray}
\label{righttail} \
&&\lim_{x \rightarrow \infty} \frac{1-\mathfrak{F}(x)}{1-\mathfrak{F}(x)
+ \mathfrak{F}(-x)} = Cp,\\
\label{lefttail}
&&\lim_{x\rightarrow \infty}\frac{\mathfrak{F}(-x)}{1-\mathfrak{F}(x)} = Cq, 
\end{eqnarray}
and $p+q=1$; the quantities $p,q$ and $C$ here are the same as in
formula (\ref{oh1}). We will say that the stable law is
\emph{unbalanced} if $p=1$ or $q=1$ above. This will happen if the
support of the variable $y$ is positive -- this will be the only case
we will consider in the sequel.

If $\chi(k)$ is the
characteristic function of our variable $y$, then the characteristic
function of the stable distribution, $\Psi(k)$, satisfies
\begin{equation}
\label{limeq} \Psi(k) = \lim_{n\rightarrow \infty}\Psi_n(k), 
\end{equation}
where 
\begin{equation}
\label{Psin}
\Psi_n= \exp(-i b_n k) \chi^n(a_n k). 
\end{equation}
\noindent\textbf{Notation.} Throughout the paper we will use the
notation $G_n$ for the distribution function of the random variable
$Y_n$ and $G$ for the corresponding stable distribution; $g_n$ for the
density of $Y_n$ and $g$ for the stable density; $\Psi_n$ for the
characteristic function of $G_n$ and $\Psi$ for the characteristic
function of the stable distribution.

\subsection{Limiting distribution of the harmonic mean $H_n$ for $\alpha=1$}
\label{limharm}
Let us go back to the example of Section \ref{prelim}, where the
random variables $x_i$ were uniformly distributed in $[0,1]$. We will
study the limiting behavior of the distribution of quantities related
to $S_n=\sum_{j=1}^n1/x_j$.

The distribution function of the variables $y_i=1/x_i$ is given by
(\ref{defF}), which implies $p=1$ and $q=0$, see formulas
(\ref{righttail}) and (\ref{lefttail}). From the behavior of the tails
of the distribution $\mathfrak{F}$ we see that $\alpha=1$, so the
norming sequence should be taken $a_n=1/n$, $b_n=\log n$. Then the
distribution $G_n$ of the variable $Y_n$ (given by equation
(\ref{XY})) converges to a stable distribution $G$. The explicit form
of the corresponding stable density, $g$, can be obtained by taking
the Fourier transform of the characteristic function $\Psi$ in formula
(\ref{expforms}) (see also \cite[Chapter XVII]{feller}). A direct
derivation of formula (\ref{Yn}) is given in Appendix \ref{explic}:
\begin{equation}
\label{Yn}
g(y)=1/(2\pi)\int_{-\infty}^\infty e^{-iky}e^{-|k|\pi/2 -ik(\log{|k|}-1+\gamma)}\,dk,
\end{equation}
where $\gamma$ is Euler's constant.

\thm{remark} {Results of this section can be easily generalized to any
density $f$ of the random variable $x$ which satisfies
$\lim_{x\to 0}f(x)>0$. For any such distribution we obtain a
stable law with exponent $\alpha=1$.}

Next, let us analyze the limiting behavior of the harmonic mean,
$H_n$. To begin we will compute the behavior of the mean of $H_n$, 
\begin{equation}
\label{mean}
\textbf{E}(H_n) \asymp \int_{-\infty}^{\infty} \frac{1}{x+\log n}
d G_n.
\end{equation}
It turns out that for this, we do not need the explicit form of the
stable distribution $G$ of $Y_n$; it is enough to use the following
information about the behavior of the tails:
\begin{equation}
\label{tails} \lim_{x \rightarrow -\infty} x G(x) = 0, \qquad
\lim_{x \rightarrow \infty}x (1-G(x)) = 1. 
\end{equation}
These equations can be obtained from (\ref{righttail}) and
(\ref{lefttail}), see XVII.5 \cite{feller}. The exact asymptotics of
the tails are computed in \cite[Chapter 2]{iblin}. 

Let us pick a large cutoff $c_n$; we take $c_n$ to tend to $\infty$,
but in such a way that $c_n = o(\log n)$, e.g. $c_n=\sqrt{\log n}$,
and rewrite equation (\ref{mean}) as
\begin{equation}
\textbf{E}(H_n) = I_1(c_n) + I_2(c_n) + I_3(c_n),
\end{equation}
where
\begin{gather}
\label{I1}
I_1(c_n) = \int_{-\infty}^{-c_n} \frac{1}{x+\log n} dG_n,\\ 
I_2(c_n) =\int_{-c_n}^{c_n} \frac{1}{x+\log n} dG_n,\\ 
\label{I3}
I_3(c_n) =\int_{c_n}^{\infty}\frac{1}{x+\log n} dG_n.
\end{gather}
We estimate these integrals separately, using equation (\ref{tails}),
the observation that $G_n(c_n) = 0$ for $c_n<1-\log n$ and the estimate on
the convergence speed of $G_n$ to $G$ as obtained in Section
\ref{sconv}. Since we are integrating over an interval of length
bounded by a constant times $\log n$, it is more than sufficient for
the speed of convergence to the stable density to be of order
$\log^2n/n$, see Theorem \ref{convthm}. Integrating by parts, we
obtain
\begin{equation}
\label{I1a}
I_1(c_n) \asymp \int_{1-\log n}^{-c_n}\frac{1}{x+\log n} dG =
\frac{G(-c_n)}{-c_n + \log n} + \int_{1-\log n}^{-c_n}\frac{G(x)}{(x+\log n)^2} d x.
\end{equation}
The first term is seen to be $o\left(\frac{1}{c_n (-c_n + \log
n)}\right)$, so given our choice of $c_n$, it is $o(1/\log n)$.  The
integral in the right hand side of equation (\ref{I1a}) is
asymptotically (in $c_n$) smaller than
\begin{equation}
\int_{1-\log n}^{-c_n} \frac{1}{-x (x + \log n)^2} d x =
O\left(\frac{1}{\log n}\right),
\end{equation}
and therefore, $I_1 = o(1/\log n).$ To show that $I_3 = o(1/\log 
n)$, we note that the integrand is dominated by $1/\log n$, while 
$\lim_{c_n \rightarrow \infty} \int_{c_n}^\infty dG_n =0$. For $I_2$ 
we have the trivial estimate (since $1/(x+\log n)$ is monotonic 
for $x > - \log n$):
\begin{equation}
\frac{G(c_n) - G(-c_n)}{\log n + c_n} \leq I_2 \leq \frac{G(c_n)
-G(-c_n)}{\log n - c_n},
\end{equation}
from which it follows that $\lim_{n\to \infty} I_2 \log n = 1$.
To summarize, we have shown
\thm{theorem} {\label{log}
For the variable $x$ uniformly distributed in $[0,1]$, $\lim_{n\to \infty}
\mathbf{E}(\log nH_n)=1.$}
\thm{remark} {\label{log1}For a general density of $x$ satisfying
$\lim_{x\to 0}f(x)>0$, we have $Y_n=1/H_n-c\log n$, and Theorem
\ref{log} generalizes to $\lim_{n\to \infty} \mathbf{E}(\log
nH_n)={\mathfrak{C}_1}=1/c.$}

In addition, we have the following weak law of large numbers for
$H_n$:
\thm{theorem}
{
\label{weaklaw}
For any $\epsilon  > 0$, $\lim_{n \rightarrow
\infty}\mathbf{P}(|H_n \log n - 1| > \epsilon) = 0.$
}

\begin{proof} Note that $$\mathbf{P}(H_n \log n -1 > \epsilon) = 
\mathbf{P}(X_n < \frac{\log n}{1+\epsilon}),$$ 
while 
$$\mathbf{P}(H_n \log n -1 < -\epsilon) = \mathbf{P}(X_n > \frac{\log
n}{1-\epsilon}).$$
Both probabilities decrease roughly as $1/(\epsilon \log n)$ using the
estimates (\ref{tails}).
\end{proof}

The above weak law indicates that if we are to hope for a limiting
distribution for $H_n$, we need to normalize it differently than by
multiplying by $\log n$. An examination of the argument above shows
that the appropriate normalization is $H_n\log^2 n  -
\log n$. Indeed, we have the following
\thm{theorem}
{
\label{limlaw}
The distributions of the random variable $H_n\log^2 n - \log n$
converges to the variable with distribution function $1-G(-x),$ where
$G$ is the limiting (stable)  distribution (of 
exponent $\alpha=1$) of variables $Y_n = X_n - \log n$.  }

\begin{proof}
The proof is quite simple. Indeed, since
$$H_n = \frac{1}{Y_n + \log n},$$
we write
\begin{eqnarray*}
\mathbf{P}(H_n \log^2 n - \log n < a) &= \mathbf{P}(\frac{\log^2
n}{Y_n + \log n} - \log n < a)\\ &= 
\mathbf{P}(\frac{ - Y_n \log n}{Y_n + \log n} < a).
\end{eqnarray*}
Since $Y_n + \log n > 0,$ we can continue:
\begin{eqnarray*}
\mathbf{P}(\frac{ - Y_n \log n}{Y_n + \log n} < a) &=
\mathbf{P}(\frac{-a \log n}{a + \log n} < Y_n)\to
\mathbf{P}(Y_n>-a)\to 1-G(-a),
\end{eqnarray*}
where we have assumed that $n$ is large enough that $a + \log n > 0.$
\end{proof}
\thm{remark}
{
\label{limlawxx} For a general density of $x$ satisfying
$\lim_{x\to 0}f(x)>0$, Theorem \ref{limlawxx} can be generalized in
the following way: the random variable $\log n(H_n\log n - {\mathfrak
C}_1)$ converges in distribution to a variable with distribution
function $1-G(-x/{\mathfrak C}_1^2),$ where $G$ is the
limiting distribution (of exponent $\alpha=1$) of variables $Y_n = X_n
- c\log n$ and ${\mathfrak C}_1=1/c$.  }

Theorem \ref{limlaw} could be viewed as a kind of an extension of Zolotarev's
identity (see \cite[Chapter XVII, Section 6]{feller} and
\cite[Theorem 2.3.4]{iblin}):
\begin{quotation}
Let $\alpha > 1.$ Then the density $p(x; \alpha)$ of the unbalanced
stable law satisfies
\begin{equation}
\label{zolotarev}
x p(x; \alpha) = x^{-\alpha} p(x^{-\alpha}; \frac{1}{\alpha}).
\end{equation}
\end{quotation}

\section{Limiting distribution of the harmonic mean $H_n$ for $\alpha\ne 1$}
\label{betan1}
Let us consider other types of the distribution of
the variable $x$ and study the limiting behavior of the corresponding
harmonic mean. If the density $f$ of the variable $x$ behaves as
\begin{equation}
\label{xbeta}
f(x)\sim x^\beta
\end{equation}
near $x=0$, then we have for the density of $y=1/x$:
$d\mathfrak{F}(y)/dy\sim |y|^{-(\beta+2)}$ as $|y|\to\infty$, which
gives $\alpha=\beta+1$ as the exponent of the stable law. Using the
material of Section \ref{stablaw} and the definition of $H_n$, we
obtain:
\begin{equation}
Y_n=\begin{cases}
    \frac{n^{1-1/\alpha}}{H_n} & 0<\alpha<1, \\
    n^{1-1/\alpha}\left(\frac{1}{H_n}-{\mathcal E}\right) & \alpha>1
  \end{cases}
\end{equation}
(here ${\mathcal E}\equiv \mathbf{E}(y)$). 

\subsection{The case $\beta > 0$}
\thm{theorem}{\label{exp}If $\beta > 0$, then  $\mathbf\lim_{n\rightarrow
 \infty} \mathbf{E} (H_n) = 1/{\mathcal E}$.}
\begin{proof} If $\beta>0$ (i.e. $\alpha>1$), then we have
$$
\lim_{n\to\infty}\mathbf{E} (H_n)=\lim_{n\to\infty}\mathbf{E}\left(\frac{Y_n}{n^{1-1/\alpha}}+{\mathcal E}\right)^{-1}=\lim_{n\to\infty}\int_{-\infty}^\infty \frac{dG_n}{\frac{x}{n^{1-1/\alpha}}+{\mathcal E}}=\frac{1}{{\mathcal E}}.
$$
\end{proof}
There is also the following Weak Law of Large Numbers:
\thm{theorem}{\label{weakal} If $\beta$ in equation (\ref{xbeta}) is
positive, then
\begin{equation*}
\lim_{n\rightarrow \infty} \mathbf{P}(\vert H_n 
-\frac{1}{{\mathcal E}}\vert> \epsilon) = 0.
\end{equation*}
}
\begin{proof} We have
\begin{eqnarray*}
\mathbf{P}\left(\vert H_n-\frac{1}{{\mathcal E}}\vert
>\epsilon\right)&=\mathbf{P}\left(\frac{1}{{\mathcal E}}-\frac{1}{{\mathcal E}+\frac{Y_n}{n^{1-1/\alpha}}}>\epsilon\right)\\
&=\mathbf{P}\left(Y_n>n^{1-1/\alpha}\frac{{\mathcal E}^2\epsilon}{1-{\mathcal E}\epsilon}\right).
\end{eqnarray*}
Since $\alpha=\beta+1>1$, then in the limit $n\to \infty$ this
quantity tends to zero.
\end{proof}

In fact, we can use a manipulation akin to that in the proof of
Theorem \ref{limlaw} to show:
\begin{theorem}
\label{limlaw2}
The random variable $n^{1-1/\alpha} (H_n - \frac{1}{{\mathcal E}})$
converges in distribution to a variable with distribution function
$1-G(-x {\mathcal E}^2)$, where the 
distribution $G$ is the unbalanced stable distribution of exponent $\alpha$.
\end{theorem}
\begin{proof}
\begin{eqnarray*}
\mathbf{P}\left(n^{1-1/\alpha} (H_n - \frac{1}{{\mathcal E}}) < a\right) &=
\mathbf{P} \left(n^{1-1/\alpha}\left(\frac{1}{{\mathcal E} + Y_n n^{1/\alpha - 1}} -\frac{1}{{\mathcal E}}\right) < a \right) \\
&= \mathbf{P}\left(\frac{Y_n}{{\mathcal E} + Y_n n^{1/\alpha
- 1}} > - a {\mathcal E}\right).
\end{eqnarray*}
The quantity ${\mathcal E} + Y_n n^{1/\alpha - 1}$ is positive because
$Y_n\ge n^{1-1/\alpha}(1-{\mathcal E})$, so we can write
\begin{eqnarray*}
\mathbf{P}\left(\frac{Y_n}{{\mathcal E} + Y_n n^{1/\alpha
- 1}} > - a {\mathcal E}\right)
&= \mathbf{P}\left(Y_n>\frac{-a{\mathcal E}^2}{1+a{\mathcal E}n^{1/\alpha
- 1}}\right)\to \mathbf{P}\left(Y_n>-a{\mathcal E}^2\right)\\
&\to 1-G(-a{\mathcal E}^2),
\end{eqnarray*}
where we have assumed that $n$ is large enough that $1+a{\mathcal
E}n^{1/\alpha - 1} > 0.$
\end{proof}

\subsection{The case $\-1 < \beta < 0$}
\label{betam}
\begin{theorem}
\label{negbeta}
 For $-1 < \beta < 0$, there
 is a constant $\mathfrak{C}_2$ such that $ \mathbf{E} (H_n/n^{1 -
 1/(\beta + 1)}) =\mathfrak{C}_2.$ 
\end{theorem}

\begin{proof}
For $-1<\beta<0$ (or $0<\alpha<1$) we would like to reason as follows:
\begin{equation}
\label{dodgy}
\lim_{n\to\infty}\mathbf{E}\left(\frac{H_n}{n^{1-1/\alpha}}\right) =\lim_{n\to\infty}\mathbf{E}\left(
\frac{1}{Y_n}\right)=\lim_{n\to\infty}\int_{-\infty}^\infty\frac{dG_n}{x}=\int_{-\infty}^\infty\frac{dG}{x}. 
\end{equation}

Since the function $1/x$ is unbounded, the weak convergence of
the distributions $G_n$ to the stable distribution $G$ is not enough
to justify the last step equality in the sequence (\ref{dodgy})
above. To justify it we need the following Lemmas:

\begin{lemma}
\label{dodlem}
Let $y_1, \dots, y_n$ be \emph{positive} independent identically distributed
random variables. Let $S_n = \sum_{i=1}^n y_i.$ Then, 
$$\mathbf{P}(S_n < a)  \leq \left[\mathbf{P}(y_1 < a)\right]^n.$$
\end{lemma}

\begin{proof}
Note that $S_n \geq \max_{1\leq i \leq n} y_i.$ 
\end{proof}

Now, in our case 
$$G_n(a) = \mathbf{P}\left(\sum_{i=1}^n y_i < a
n^{\frac{1}{\alpha}}\right) \leq
\left[\mathbf{P}\left(y_1 < a n^{\frac{1}{\alpha}}\right)\right]^n =
\left[\mathbf{P}\left(x_1 > \frac{1}{a n^{\frac{1}{\alpha}}}\right)\right]^n,
$$
where the inequality follows from Lemma \ref{dodlem} (and recall that
$x_i = 1/y_i$).

The probability $\mathbf{P}(x_1 > b)$ has the following
properties:
\begin{itemize}
\item[A)] $\mathbf{P}(x_1 > b) = 0$ for $b \geq 1,$
\item[B)] 
$1-\mathbf{P}(x_1 > b) \sim c b^\alpha,$ for $b \ll 1,$
\item[C)]
$\mathbf{P}(x_1 > b) < 1$ for $b > 0.$
\end{itemize}
\begin{lemma}
\label{a1}
$G_n(a) = 0$ for $a \leq n^{-\frac{1}{\alpha}}.$
\end{lemma}
\begin{proof}
Follows from the definition of $G_n$ and Property A.
\end{proof}
\begin{lemma}
\label{a2}
There exists a $b_0$ such that $1-\mathbf{P}(x_1 > b) < 2 c' b^\alpha$
for all $b < b_0$ for some $c^\prime > 0.$
\end{lemma}
\begin{proof}
This follows from Property B, with $c^\prime = 2 c.$
\end{proof}
\begin{lemma}
\label{a3}
If $a n^{\frac{1}{\alpha}} > 1/b_0$ ($b_0$ as in the statement of
Lemma \ref{a2}), then
$$G_n(a) \leq (1 - c^\prime a^{-1/\alpha} n^{-1})^n \sim
\exp(-c^\prime a^{-1/\alpha}).$$
\end{lemma}
\begin{proof}
Follows from Lemma \ref{dodlem}.
\end{proof}
\begin{lemma}
\label{a4}
$G(1/(b_0 n^{1/\alpha})) \leq [\mathbf{P}(x_1 > b_0)]^n.$ 
\end{lemma}
\begin{proof}
Follows immediately from Lemma \ref{dodlem}.
\end{proof}

Now we write: 
\begin{eqnarray*}\int_0^\infty \frac{d G_n}{x} &= 
\left(\int_0^{n^{-\frac{1}{\alpha}}} + 
\int_{n^{-\frac{1}{\alpha}}}^{n^{-\frac{1}{\alpha}}/b_0} + 
\int_{n^{-\frac{1}{\alpha}}/b_0}^{C} + \int_{C}^\infty 
\right) \frac{d G_n}{x}\\ &= 
I_0(n) + I_1(n) + I_2(n) + I_3(n).
\end{eqnarray*} 
To analyze the above decomposition, we should first belabor the
obvious:
\begin{lemma}
\label{b0}
$$
\int_a^b \frac{dG_n}{x} = \frac{G_n(b)}{b} - \frac{G_n(a)}{a} + 
\int_a^b \frac{G_n}{x^2} d x.
$$
\end{lemma}
\begin{proof}
Integration by parts.
\end{proof}
\begin{lemma}
\label{b1}
$I_0(n) = 0.$
\end{lemma}
\begin{proof}
The integrand vanishes in the interval by Lemma \ref{a1}.
\end{proof}
\begin{lemma}
\label{b2}
$$\lim_{n \rightarrow \infty} I_1(n)= 0.$$
\end{lemma}
\begin{proof}
By Lemma \ref{a4}, $G_n < [\mathbf{P}(x_1 > b_0)]^n.$ The result follows by integration
by parts (Lemma \ref{b0}).
\end{proof}
\begin{lemma}
\label{b3}
$$\lim_{n \rightarrow \infty} I_2(n) \leq
\frac{\exp\left(C^{-\frac{1}{\alpha}}\right)}{C} +
\int_0^{C} \frac{\exp\left(x^{-\frac{1}{\alpha}}\right)}{x^2} d x.
$$
\end{lemma}
\begin{proof}
Follows from Lemma \ref{b0} and Lemma \ref{a3}.
\end{proof}
\begin{lemma}
\label{b4}
$$\lim_{n \rightarrow \infty} I_3(n) = \int_{C}^\infty \frac{d
G}{x}.$$
\end{lemma}
\begin{proof}
This follows from the weak convergence of $G_n$ to $G.$
\end{proof}
The derivation (\ref{dodgy}) is justified. Indeed, if we make the
constant $C$ above large, we see that the integral of $d G_n/x$ is
bounded, hence so is the integral of $d G/x.$ Convergence follows from
the dominated convergence theorem (or by making $C$ small). We have
incidentally shown that the density of the stable law decays
exponentially as $x \rightarrow 0^+$ (exact expression can be found in
\cite[Chapter 2]{iblin}),
\end{proof}
\begin{theorem}
\label{limlaw3}
The quantity $H_n/n^{1-1/\alpha}$
converges in distribution to the variable with distribution function
$1-G(1/x)$, where $G$ is the unbalanced stable law of exponent $\alpha.$
\end{theorem}
\begin{proof}
The proof is immediate.
\end{proof}

\section{\label{classpol}A class of random polynomials}
Let $x_1,\ldots,x_n$ be independent identically distributed random
variables with values between zero and one. Let us consider
polynomials whose roots are located at $x_1,\ldots,x_n$:
\begin{equation}
\label{poly}
p(x)= \prod_i^n(x-x_i)=x^n + \sum_{i=0}^{n-1} c_i x^i.
\end{equation}
Given the distribution of $x_i$, we would like to know the
distribution law of the roots of the derivatives of $p(x)$.
\subsection{Uniformly distributed roots}
Let us denote the roots of $\frac{dp(x)}{dx}\equiv p'(x)$ by $\mu_i$, $1\le
i\le n-1$, and assume that $\mu_i\le \mu_{i+1}$ for all $i$. It is
convenient to denote the smallest of $x_j$ by $m_1$, i.e. $m_1\equiv
\min_j x_j$, the second smallest of $x_j$ as $m_2$ and so forth. It is
clear that
\begin{equation}
\label{m1mum2}
m_i\le \mu_i\le m_{i+1},\,\,\,\,\,\,\,1\le i\le n-1.
\end{equation}

We now assume that the $x_j$ are independently uniformly distributed in
$[0, 1]$. The distribution of $m_1$ is easy to compute: the
probability that $m_1> \alpha$ is simply the probability that all of
the $x_j$ are greater than $\alpha$, which is to say,
\begin{equation*}
P(m_1>\alpha)=(1-\alpha)^n.
\end{equation*}
Using this distribution function, one can show that
\begin{equation*}
{\mathbf E} (m_1)=\frac{1}{n+1}.
\end{equation*}
In fact, it is not hard to see that ${\mathbf E}(m_i) = i/(n+1);$ the
reader may wish to consult \cite{feller} (page 34). We thus have:
\begin{equation}
\label{mui1}
\frac{i}{n+1}\le {\mathbf E}( \mu_i) \le \frac{i+1}{n+1},\,\,\,\,1\le i\le n-1.
\end{equation}
In particular, for large values of $n$ we have the estimate
\begin{equation*}
{\mathbf E}( \mu_*) \sim \frac{1}{n},
\end{equation*}
where the notation $\mu_*$ is used for the {\it smallest} root of the
derivative.
\subsection{\label{rootzero} More precise locations of roots of the 
derivative, given that the smallest root of the polynomial is fixed}
In the previous section we have noted that if the roots of $p (x)$ are
distributed uniformly in $[0, 1]$, then so are the roots of
$p^\prime(x)$. In order to understand better the
distribution of the roots of $p^\prime(x)$, first let
$$
p(x) = (x-x_1) \dots (x - x_n),
$$
then we can write
$$
p^\prime(x) = p(x) \sum_{j=1}^n \frac{1}{x-x_j}.
$$
In the generic case where $p(x)$ has no multiple roots, a root
$\mu$ of $p^\prime(x)$ satisfies the equation
\begin{equation}
\label{drooteq} \sum_{j=1}^n \frac{1}{x_j - \mu} = 0.
\end{equation}
This was interpreted by Gauss (in the more general context of complex
roots) as saying that $\mu$ is in equilibrium in a force field where
force is proportional to the inverse of distance, and the ``masses''
are at the points $x_1, \dots, x_n.$ Gauss used this simple
observation to deduce the Gauss-Lucas theorem to the effect that the
zeros of the derivative lie in the convex hull of the zeros of the
polynomial (see \cite{marden}). We will use it to get more precise
location information on the zeros.  In particular, consider the
smallest root $\mu_*$ of $p^\prime(x)$. It is attracted from the left
only by the root $x_1$ of $p$, and from the right by all the other
roots, so we see

\thm{lemma}{\label{trivest}  For all $2\le
i\le n$, $(m_1 + m_i)/2 \geq \mu_*,$ with equality if and only if
$n=2.$ }

\thm{remark} {In the sequel, we shall assume that the smallest root of
$p$ equals zero (i.e. $m_1=0$; for simplicity of notation we assume
that $x_n=0$).}

Inequalities (\ref{m1mum2}) still give a good estimate for the roots
$\mu_2,\ldots,\mu_{n-1}$. One can similarly show that
\begin{equation}
\label{mui2}
\frac{i-1}{n}\le {\mathbf E}( \mu_i) \le \frac{i}{n},\,\,\,\,2\le i\le n-1.
\end{equation}
However, for $\mu_1=\mu_*$, inequalities (\ref{m1mum2}) give $0\le
\mu_*\le m_2$. For uniformly distributed $x_i$ this tells us that
$\mu_*$ decays like $1/n$ or faster. We would like to obtain a more
precise estimate for the large $n$ behavior of $\mu_*$.

The random polynomial, $p(x)$, now has the form
\begin{equation*}
p(x)=x\prod_{i=0}^{n-1}(x-x_i)=x\sum_{i=0}^{n-1}c_ix^i.
\end{equation*}
We need to estimate the smallest root of $p'$, $\mu_*$.

\begin{theorem}
\label{bounds}
\begin{equation*}
\frac{1}{2} \left(\sum_{i=1}^{n-1}1/x_i\right)^{-1}\le
\mu_*\le\left(\sum_{i=1}^{n-1}1/x_i\right)^{-1}.
\end{equation*}
\end{theorem}
%
\begin{proof} 
The smallest root $\mu_*$ satisfies the equation 
\begin{equation}
\label{mueq} 1/\mu_*= \sum_{i=1}^{n-1} 1/(x_i - 
\mu_*). 
\end{equation}
By Lemma \ref{trivest},
\begin{equation}
\label{stupideq}
 \frac{x_i}{2} \leq x_i - \mu_* \le x_i,\quad\quad 1\le i\le n-1.
\end{equation} 
The result follows immediately from equations (\ref{mueq}) and
(\ref{stupideq}). \end{proof}

Now it is clear that in order to find an estimate for $\mu_*$, we need
to study the behavior of $\left(\sum_{i=1}^{n-1}1/x_i\right)^{-1}$. In
terms of the harmonic mean of independent random variables
$x_1,\ldots, x_{n-1}$, we have
\begin{equation}
\label{relate}
\frac{1}{2(n-1)}H_{n-1}\le \mu_*\le \frac{1}{n-1}H_{n-1},
\end{equation}
or, for large values of $n$,
\begin{equation}
\label{mustar}
{\mathbf E}( \mu_*) \sim \frac{1}{n} {\mathbf E}( H_n).
\end{equation}
Using Theorems \ref{log} and \ref{weaklaw}, we readily obtain
\begin{equation*}
{\mathbf E}( \mu_*) \sim \frac{1}{n\log{n}}
\end{equation*}
\section{\label{classmat}A class of stochastic matrices}
\label{stochmat}
Let $T$ be an $n$ by $n$ matrix constructed as follows:
\begin{equation}
\label{stoch}
T_{ij}=\left\{ \begin{array}{l}
(1-a_{i})/(n-1),\,\,\,\,\,\,i\ne j,\\
a_{i},\hspace{3.0cm}i=j,
\end{array}\right.
\end{equation}
where the numbers $a_i$ are independently distributed between $0$
and $1$. We want to study the large $n$ behavior of the second
largest eigenvalue of $T$ (the largest eigenvalue is equal to 1).
We will denote this eigenvalue as $\lambda_*$. In the next section
we will provide some motivation for this choice of stochastic
matrices.

\subsection{Motivation: the memoryless learner algorithm}
\label{memoryless}
The following is a typical learning theory setup (see \cite{niy}). We
have $n$ sets (which we can think of as \textit{concepts}),
$R_1,\ldots, R_n$. Each set $R_i$ is equipped with a probability
measure $\nu_i$ The \textit{similarity matrix} $A$ is defined by
$a_{ij}=\nu_i(R_j)$. Since the $\nu_i$ are \textit{probability}
measures, we see that $0\le a_{ij}\le 1$ and $a_{ii}=1$ for all
$i,j$. Now, the teacher generates a sequence of $N$ examples referring
to a single concept $R_k$, and the task of the student is to guess the
$k$ (i.e. to {\it learn} the concept $G_k$), hopefully with high
confidence.

The learner has a number of algorithms available to him. For instance,
the student may decide in advance that the concept being explained is
$R_1$, and ignore the teacher's input, insisting forever more that the
concept is $R_1$. While this algorithm occasionally results in
spectacular success, the probability of this is independent of the
number of examples, and is inversely proportional to the number $n$ of
available concepts. Here we will consider a more practical and
mathematically interesting algorithm, namely,

\medskip\noindent \textbf{Memoryless learner algorithm.}  The student
picks his initial guess at random, as before. However, now he
evaluates the teacher's examples, and if the current guess is
incorrect (i.e. if the teacher's example is inconsistent with the
current guess), he switches his guess at random. The name of the
algorithm stems from the fact that the student keeps no memory of the
history of his guesses, and will occasionally switch his guess to one
previously rejected.

It is clear that with the memoryless learner algorithm, the student
will never be able to learn the set $R_k$ if $R_k \subset R_l$. We
call such a situation \textit{unlearnable}, and do not consider it in
the sequel. In terms of the similarity matrix, this can be rephrased
as the assumption that $a_{ij} < 1, \quad i \neq j.$

To define our mathematical model further, we will assume that the
student picks the initial guess uniformly: ${\bf
p}^{(0)}=(1/n,\ldots,1/n)^T.$ The discrete time
evolution of the vector ${\bf p}^{(t)}$ is a Markov process with
transition matrix $T^{(k)}$, which depends on the teacher's concept,
$R_k$, and the similarity matrix, $A$. That is:
\begin{equation}
\label{Tij}
T^{(k)}_{ij}=\left\{ \begin{array}{l}
(1-a_{ki})/(n-1),\,\,\,\,\,\,i\ne j,\\
a_{ki},\hspace{3.0cm}i=j.
\end{array}\right.
\end{equation}
After $N$ examples, the probability that the student believes that the
correct concept is $R_j$ is given by the $j$th component of the vector
$({\bf p}^{(N)})^T=({\bf p}^{(0)})^T(T^{(k)})^N$. In particular, the
probability that the student's belief corresponds to reality (that is,
$j = k$) is given by:
\begin{equation}
\label{Qij}
Q_{kk}(N)=[({\bf p}^{(0)})^T\,(T^{(k)})^N]_k.
\end{equation}
It is clear that the dynamics of the memoryless learner algorithm 
is completely encoded by the  matrix $T$ defined
above by (\ref{Tij}).

We are interested in the rate of convergence as a function of $n$, the
number of possible concepts. We define the convergence rate of the
algorithm as the rate of the convergence to $0$ of the difference
\begin{equation*}
1-Q_{kk}(N).
\end{equation*}
In order to simplify notation, let us set $k=1$ and skip the
corresponding subscript/superscript. In order to evaluate the
convergence rate of the memoryless learner algorithm, let us represent
the matrix $T^{(1)}\equiv T$ as follows:
\begin{equation}
\label{oh}
T=V\Lambda W,
\end{equation}
where the diagonal matrix $\Lambda$ consists of the eigenvalues of
$T$, which we call $\lambda_i$, $1\le i\le n$; representation
(\ref{oh}) is possible in the generic case. The columns of the matrix
$V$ are the right eigenvectors of $T$, ${\bf v}_i$. The rows of the
matrix $W$ are the left eigenvectors of $T$, ${\bf w}_i$, normalized
to satisfy $<{\bf w}_i, {\bf v}_j>=\delta_{ij}$, where $\delta_{ij}$
is the Kronecker symbol (so that $VW=WV=I$). The eigenvalues of $T$
satisfy $|\lambda_i|\le 1$. We have
\begin{equation*}
T^N=V\Lambda^NW.
\end{equation*}
Let us arrange the eigenvalues in decreasing order, so that
$\lambda_1=1$ and $\lambda_2\equiv \lambda_*$ is the second largest
eigenvalue (we assume that it has multiplicity one). If $N$ is large,
we have $\lambda_i^N\ll \lambda_*^N$ for all $i\ge 3$, so only the
first two largest eigenvalues need to be taken into account. This
means that in order to evaluate $T^N$ we only need the following
eigenvectors: ${\bf v}_1=(1/n,1/n,\ldots, 1/n)^T$, ${\bf v}_2$, $ {\bf
w}_1=(n,0,\ldots,0)$, and $ {\bf w}_2$ (it is possible to check that
the contribution from the other components contains multipliers
$\lambda_i^N$ with $i>2$ and thus can be neglected, see the
computation for $C_n$ in Section \ref{eig}). It follows that
\begin{equation}
\label{conv}
Q_{11}=1-C_n\lambda_*^N,
\end{equation}
where
\begin{equation}
\label{constC}
C_n=-\frac{1}{n}\sum_{j=1}^n[{\bf v_2}]_j[{\bf w}_2]_1.
\end{equation}
The convergence rate of the memoryless
learner algorithm can be found by estimating $\lambda_*$ and
$C_n$. It turns out that a good understanding of  $\lambda_*$ (Section
\ref{mustar1}), helps us also estimate $C_n$ (this is done in Section
\ref{eig}).

\subsection{\label{mustar1}Second largest eigenvalue and the smallest
root of the derivative of a random polynomial} 
\label{charderiv}
Let $Z = I - T,$ and let $x_i = 1- a_i.$ The matrix $Z$ satisfies
$Z_{ii} = x_i$, while $Z_{ij} = -\frac{x_i}{n-1},$ for $j\neq i$. We have 
\begin{equation}
\label{Zdef}
Z = \frac{n}{n-1} D_x(I - \frac{1}{n} J_n), 
\end{equation}
where $J_n$ is the $n \times n$ matrix of all ones, and $D_x$ is a
diagonal matrix whose $i$-th element is $x_i$. It is convenient to
introduce the matrices
\begin{equation}
\label{Mn}
M_n=I-\frac{1}{n} J_n=\left(\begin{array}{cccc}
1 & -1/n & -1/n & \ldots \\
-1/n & 1 & -1/n & \ldots \\
-1/n & -1/n & 1 & \ldots \\
\ldots & & & 
\end{array}\right)
\end{equation}
and 
\begin{equation}
\label{Bdef}
B = D_x M_n.
\end{equation}
The second largest eigenvalue of $T$, which we denote as $\lambda_*$,
and the smallest nontrivial eigenvalue of $B$, $\mu_*$, are related as
\begin{equation}
\label{connect}
\lambda_*=1-\frac{n}{n-1}\mu_*.
\end{equation}

In what follows we will write down the characteristic polynomial of
$B$. Let us recall the following 
\begin{fact}
\label{minorfact}
Let $A$ be an $n\times n$ matrix. Then the coefficient of $x^{n-k}$ in
the characteristic polynomial $p_A(x)$ of $A$ (defined to be $\det(x I_n
-A)$) is given by 
{\rm $$\sum_{\tiny \begin{array}{c}\mbox{$k$-element subsets}\\\mbox{$S$ of $\{1, \dots, n\}$}\end{array}}
(-1)^k\det m_{\bar{S}},$$}
where $m_{\bar{S}}$ is the matrix obtained from $A$ by
deleting those rows and columns whose indices are \emph{not} elements
of $S$ (we call $m_S$ the minor of $A$ corresponding to $S$).
\end{fact}

We will need the following lemmas:

\begin{lemma}
\label{minorlem}
Let $A$ be an $n \times n$ matrix, and let $D_x$ be as above. Let
$m_{i_1, \dots, i_k}$ be the $i_1, \dots, i_k$ minor of $M$ (that
is, the sub-matrix of $i_1, \dots, i_k$ rows and
$i_1, \dots, i_k$ columns of $M$). Then the minor $\gamma_{i_1, \dots,
i_k}$ of the matrix $C = D_x A$ satisfies 
$$\det \gamma_{i_1, \dots,i_k} = \det m_{i_1, \dots,i_k} \prod_{l=1}^k
x_{i_l}.$$
\end{lemma}

\begin{proof}
This is immediate, since the $j$-th row of $C$ is $x_j$ times the
$j$-th row of $A$.
\end{proof}

\begin{lemma}
\label{charpolylem}
The characteristic polynomial of $M_n = I - \frac{1}{n} J_n$ equals 
$x(1-x)^{n-1}.$
\end{lemma}
\begin{proof}
Immediate, since the bottom eigenvalue of $M_n$ (given in equation
(\ref{Mn})) is zero and the rest are $1$.
\end{proof}
\begin{lemma}
\label{minorlem2}
All the $k\times k$ minors of $M_n$ (defined in the statement of
Lemma \ref{charpolylem}) are equal.
\end{lemma}
\begin{proof}
By inspection -- all the $k \times k$ minors of $M_n$ are $M_k +
\left(\frac{1}{n} - \frac{1}{k}\right)J_k$.
\end{proof}
\begin{lemma}
\label{detform}
The determinants $d_k$of the $k \times k$ minors of $M_n$ are equal to
$\frac{k}{n}.$
\end{lemma}
\begin{proof}
We know that the $\binom{n}{k} d_k = \binom{n-1}{k-1},$ from Lemmas
\ref{charpolylem} and \ref{minorlem2}. From this the assertion follows
immediately. 
\end{proof}
Now, let 
$$p_D(x) = x^n + \sum_{i=0}^{n-1} c_i x^i$$
be the characteristic polynomial of $D_x.$

\thm{lemma}{The characteristic polynomial of $B$, $p_B(x)$, is given
by:
\begin{equation}
\label{charpoly}
p_B(x) = x^n + \sum_{i=0}^{n-1} \frac{i}{n} c_i x^i,
\end{equation}
where $c_i$ are as above.}

\begin{proof} 
>From Lemma \ref{minorlem} combined with Lemma \ref{detform}, we see
that the coefficient of $x^i$ in $p_B(x)$ is given by 
$$\frac{i}{n}\sum_{\tiny \begin{array}{c}\mbox{$i$-element subsets}\\\mbox{$S$ of $\{1, \dots, n\}$}\end{array}}
\prod_{j\in S} x_j.$$

The sum is just the $i$-th elementary symmetric function of the $x_1,
\dots, x_n,$ which is equal to $c_i$. The assertion follows.
\end{proof}

Notice that the constant term of $p_B$ vanishes, so we can write
\begin{equation*}
p_B(x) = x q(x),
\end{equation*}
where
\begin{equation*}
q(x)=x^{n-1} + \sum_{i=0}^{n-2} \frac{i+1}{n} x^i.
\end{equation*}
But obviously
\begin{equation*}
q(x) = \frac{1}{n}\frac{d p_D(x)}{d x}, 
\end{equation*}
so we have 
\thm{lemma}{\label{amuz} The characteristic polynomials of the matrix
$B$ defined in (\ref{Bdef}) and the diagonal matrix $D_x$ with
elements $x_i$ are related by 
$$ p_B(x) = \frac{x}{n} p_D^\prime(x).$$
}

This relates the eigenvalues of the matrix $B$ and the zeros of the
polynomial $q(x)$ (and $p_D^\prime(x)$). In its turn, the smallest
eigenvalue of $B$ is related to the second largest eigenvalue of our
matrix $T$ by equation (\ref{connect}).

We can see that studying the second largest eigenvalue of a stochastic
matrix of class (\ref{stoch}) is reduced to the problem of the
smallest root of the derivative of the stochastic polynomial of class
(\ref{poly}), with $x_i=1-a_i$. Note that by the definition of matrix
$T^{(k)}$, one of the quantities $1-a_{ki}=x_i$ is equal to zero. This
means that in order to find the distribution of the second largest
eigenvalue of such a matrix, we need to refer to Section
\ref{rootzero}, i.e. the case where one of the roots of the random
polynomial was fixed to zero, and the rest were distributed uniformly.
\subsection{\label{eig}Eigenvectors of stochastic matrices}
Next, let us study eigenvectors of stochastic matrices, and derive an
estimate for $C_n$ in equation (\ref{constC}). Consider the matrix $Z$
defined in equation (\ref{Zdef}). We can write $Z = W^t D_\mu V$,
where $V$ and $W$ are the matrices of right and left eigenvectors
(respectively) of $Z$, and $D_\mu$ is a diagonal matrix whose entries
are the eigenvalues of $Z$. We know that the right eigenvector of $Z$
corresponding to the eigenvalue $0$ is the vector ${\bf v}_1 = (1,
\dots, 1)^T$, while the left eigenvector is the vector ${\bf w}_1 =
(1, 0, \dots, 0).$ To write down the eigenvector ${\bf v}_i$ ($i>1$)
we write ${\bf v}_i = {\bf v}_1 + {\bf u}_i$, where $<{\bf v}_1, {\bf
u}_i> = 0$ -- we can always normalize ${\bf v}_i$ so that this is
possible. If the corresponding eigenvalue is $\mu_i$, we write the
eigenvalue equation:
\begin{equation*}
\mu_i {\bf v}_i = Z {\bf v}_i = \frac{n}{n-1} D_x (I - \frac{1}{n} J_n) ({\bf v}_1 + {\bf u}_i)
= \frac{n}{n-1} D_x {\bf u}_i.
\end{equation*}
This results in the following equations for $u_{ij}$ -- the $j$-th
coordinate of ${\bf u}_i$ (for $j > 1$; $u_{i1} = -1$):

\begin{equation}
\label{veq}
\mu_i + \mu_i u_{ij} = \frac{n}{n-1} x_j u_{ij},
\end{equation}
and so
\begin{equation}
\label{veq2}
u_{ij} = \frac{\mu_i}{\frac{n}{n-1} x_j - \mu_i}.
\end{equation}

On the other hand, the eigenvalue equation for ${\bf w}_i$ is
\begin{equation}
\label{weq}
\mu_i {\bf w}_i = Z^t {\bf w}_i = \frac{n}{n-1} (I - \frac{1}{n} J_n) D_x {\bf w}_i,
\end{equation}
resulting in the following equations for the coordinates:
\begin{equation}
\label{weq2}
\mu_i w_{ij} = \frac{n}{n-1}\left(x_j w_{ij} - \frac{1}{n} \sum_{k = 1}^n
x_k w_{ik}\right).
\end{equation}

If we assume that $x_1 = 0$, then setting $j = 1,$ we get
\begin{equation}
\mu_i w_{i1} = - \frac{n}{n-1} \frac{1}{n} \sum_{k = 1}^n x_k w_{ik},
\end{equation}
and so
\begin{equation}
\frac{1}{n} \sum_{k = 1}^n x_k w_{ik} = - \frac{n}{n-1} \mu_i w_{i1},
\end{equation}
and equation (\ref{weq2}) can be rewritten (for $j > 1$) as
\begin{equation}
\mu_i w_{ij} = \frac{n}{n-1}\left(x_j w_{ij} + \frac{n}{n-1} \mu_i
w_{i1}\right),
\end{equation}
to get
\begin{equation}
\label{weq3}
w_{ij} = \frac{\frac{n}{n-1} \mu_i w_{i1}}{- \frac{n}{n-1}x_j
 + \mu_i}
\end{equation}
Now, let us assume that $i = 2$, and in addition $\mu_2 \ll x_k, \quad
k>1$. While it follows immediately from Lemma \ref{trivest} that
$\mu_2 < x_2$, we comment that by our Weak Law of Large Numbers
(Theorem \ref{weaklaw}), the probability that $\mu_2 > c n/\log n$
goes to 0 with $c$, whereas the probability that $|x_k - (k-1)/n| >
c_2/n$ goes to zero with $c_2$ (detailed results on the distribution
of order statistics can be found in \cite[Chapter I]{feller}).

\thm{remark}{ The assumption that $\mu_* \ll x_2$ is least
justified if we have reason to believe that $x_2 \ll 1/n$. }
Thus we can write approximately:
\begin{equation}
\label{veq2a}
v_{2j} \approx \frac{\mu_2}{ x_j} + 1,
\end{equation}
while
\begin{equation}
\label{weq3a}
w_{2j} \approx - \frac{\mu_2 w_{21}}{x_j}
\end{equation}
Since we must have $<w_2, v_2> = 1$, we have:
\begin{equation*}
- w_{21} \sum_{j=2}^n \frac{\mu_2}{x_j}(\frac{\mu_2}{x_j} + 1) =
1, \end{equation*}
which implies that
\begin{equation*}
w_{21} \approx -\frac{1}{\mu_2} \frac{1}{\sum_{j=2}^{n}\frac{1}{x_j}},
\end{equation*}
which, in turn, implies that $1/2 \leq |w_{21}| \leq 1.$ This means
that we have the following estimate for the quantity $C_n$ in
(\ref{constC}):
\begin{equation*}
1 \leq -\frac{1}{n}w_{21} \sum_{j=2}^n v_{2j} \leq 2, 
\end{equation*}
i.e. 
\begin{equation}
\label{consC2}
1\le C_n\le 2.
\end{equation}
\subsection{\label{convmem}Convergence of the memoryless learner algorithm}
Let us assume that the overlaps between concepts, $a_{ki}\equiv a_i$
in the matrix $T$, are independent random variables distributed
according to density $\tilde{f}(a)$. Then the variables $x_i=1-a_i$
have the probability density $f(x)=\tilde{f}(1-x)$. Our results for
the rate of convergence of the memoryless learner algorithm can be
summarized in the following
\thm{theorem}{\label{main} Let us assume that the density of overlaps,
$f(x)$, approaches a nonzero constant as $x\to 0$. Then in order for
the learner to pick up the correct set with probability $1-\Delta$, we
need to have at least
\begin{equation}
\label{last}
N_\Delta\sim |\log \Delta|(n\log n). 
\end{equation}
sampling events.}
\begin{proof} Combining equations (\ref{conv}), (\ref{consC2}) and (\ref{connect}) we
can see that in order for the learner to pick up the correct set with
probability $1-\Delta$, we need to have at least
\begin{equation}
\label{Nandmu}
N_\Delta\sim|\log \Delta|/\mu_*
\end{equation}
sampling events. Since $\beta=0$ (see equation (\ref{xbeta}), we have
$\alpha=1$. Using bounds (\ref{relate}) which relate $\mu_*$ to the
harmonic mean, and the weak law of large numbers (Theorem
(\ref{weaklaw})), we obtain estimate (\ref{last}).  This estimate
should be understood in the following sense: as $n\to\infty$, the
probability that the ratio $\mu_*^{-1}/(n\log n)$ deviates from $1$ by
a constant amount, tends to zero. Therefore, the right hand side of
(\ref{Nandmu}) behaves like the right hand side of (\ref{last}) with
probability which tends to one as $n$ tends to infinity.
\end{proof}

For other distributions we have 
\thm{theorem}{\label{main1} If the probability density of overlaps,
$f(x)$, is asymptotic to $x^\beta + O((x)^{\beta +\delta}), \quad
\delta, \beta > 0$, as $x$ approaches $0$, then 
$$N_\Delta \sim |\log\Delta|n;$$
if $-1 < \beta < 0$, then 
$$\lim_{x \rightarrow \infty}  \mathbf{P}\left(\frac{1}{x} < \frac{N_\Delta}{|\log\Delta|n^{1/(1+\beta)}} < x\right) = 1.$$
}
\begin{proof} The proof uses the results on the harmonic mean in
section \ref{betan1}.
\end{proof}

\bibliographystyle{alpha}

\begin{thebibliography}{xxxxxxxxxxxx}

\bibitem[BenOrsz]{benorsz}
C.~M.~Bender and S.~Orszag (1999) \textit{Advanced mathematical
methods for scientists and engineers, I,\/} Springer-Verlag, New
York.

\bibitem[FellerV2]{feller}
Feller,W. (1971)
\textit{An Introduction to Probability Theory and its Applications,\/}
vol. 2,  John Wiley and Sons.

\bibitem[KNN2001]{knn} 
Komarova, N.~L., Niyogi,~P. and
Nowak,~M.~A. (2001) The evolutionary dynamics of grammar acquisition,
\textit{J.~Theor.~Biology}, {\bf 209}(1), pp. 43-59.

\bibitem[KN2001]{kn} 
Komarova, N.~L. and Nowak, M.~A. (2001) Natural
selection of the critical period for grammar acquisition, {\it
Proc. Royal Soc. B}, {\bf 268}, pp. 1189-1196.

\bibitem[Hall1981]{phall}
Hall,~P. (1981). On the rate of convergence to a stable law. 
\textit{Journal of the London Mathematical Society}(2), 
\textbf{23}, pp.~179-192.

\bibitem[IbLin1971]{iblin}
Ibragimov and Linnik.
\textit{Independent and Stationary sequences of Random Variables},
Wolters-Noordhoff Publishing, Groningen, 1971. 

\bibitem[KK2001]{kus}
Kuske, R. \& Keller,J.B. (2001). Rate of Convergence to a Stable 
Law. {\it SIAM Jour. Appl. Math.}, {\bf 61}, N.4, pp. 1308-1323.

\bibitem[Marden1966]{marden}
Marden, M. (1966) \textit{Geometry of Polynomials, second
edition}, American Math.~Soc., Providence, RI.

\bibitem[Niyogi1998]{niy}
Niyogi, P. (1998). {\it The Informational Complexity of Learning}. Boston: Kluwer.

\bibitem[Norris1940]{nor}
Norris, N. (1940) The standard errors of the geometric and harmonic
means and their applications to index numbers,
\emph{Ann.~Math.~Statistics} \textbf{11}, pp 445-448.

\bibitem[NKN2001]{nkn}
Nowak, M.~A., Komarova,~N.~L., Niyogi,~P. (2001) Evolution of
universal grammar, \textit{Science} \textbf{291}, 114-118.

\end{thebibliography}

\appendix

\section{\label{explic}Derivation of the stable law for the uniform distribution of $x_i$} 

Here we will provide an explicit derivation of the stable law,
equation (\ref{Yn}), in the case where the random variables $x_i$ are
uniformly distributed between zero and one. The characteristic
function corresponding to the distribution of $y_i=1/x_i$,
equation (\ref{defF}), is given by
\begin{equation}
\chi(k) = \int_1^\infty \frac{e^{i k y}}{y^2} d y.
\end{equation}
This can be evaluated explicitly; here we only present the
computations for positive $k$, to avoid clutter. We have
\begin{equation}
\label{concr} \chi(k)=-|k|\pi/2+\cos k+k{\rm Si}\,(k)+i[\sin k-k{\rm 
Ci}\,(|k|)],
\end{equation}
where $\mathrm{Si}$ and $\mathrm{Ci}$ are sin- and cos- integrals,
respectively (this expression is obtained with \textit{Mathematica}).
The behavior of $\chi(k)$ at $0$ can be easily computed using the
above formula:
\begin{equation}
\label{concrasymp} \chi(k) = 1 - i (\gamma - 1) k - \frac{\pi}{2} k 
- i k \log k + \frac{1}{2} k^2 + o(k^2),
\end{equation}
where $\gamma\approx 0.577216$ is Euler's constant.

We can also obtain the asymptotics for $\chi(k)$ directly, as
follows. First, we change variables, and set $u = k y,$ to obtain
\begin{equation*}
\chi(k) = k \int_k^\infty \frac{e^{i u}}{u^2} d u.
\end{equation*}
Let 
$$I(k) = \int_k^\infty \frac{e^{i u}}{u^2} d u\quad\quad{\rm and}\quad \quad  I_R(k) =
\int_k^R \frac{e^{i u}}{u^2} d u.$$
Clearly, $I(k) = \lim_{R \rightarrow \infty} I_R(k).$ Since the
integrand (call it $f(u)$) has no poles, except at $0$, we see that
$I_R(k) + J(k) + L(R) - K_R(k) = 0$, where $J(k)$ is the integral of
$f(u)$ along the positive quadrant of the circle $|z| = k$, $L(R)$ the
integral along the positive quadrant of he circle $|z| = R,$ and
$K_R(k)$ is the integral along the imaginary axis, from $i k$ to $iR$.
It is easy to see that $\lim_{R \rightarrow \infty} L(R) = 0,$ while
$$k J(k) = \frac{1}{2} \int_0^{\pi/2} \exp(i k e^{i \theta} - i 2
\theta) d \theta,$$ 
the Taylor series of which is easily evaluated by
expanding the integrand in a Taylor series.  The integral $K_R(k)$
reduces to the integral of $\exp(-u)/u^2$, the asymptotics of which
can be easily obtain by repeated integration by parts.

The characteristic function of $Y_n=a_n S_n - b_n$, $\Psi_n$, can be
obtained by setting $a_n = 1/n,$ and $b_n = \log n$, and using
equation (\ref{Psin}):
\begin{equation}
\label{exampcalc}
\begin{split}
\Psi_n(k) & = \exp(- i k\log n) (\chi(k/n))^n\\ &= \exp\left( -
  \left(\frac{\pi}{2} +i (\gamma - 1)\right)k + i k\log
  k\right)\left(1 + O\left(\frac{k^2\log^2 n}{n}\right)\right),
  \end{split} \end{equation}
where we expanded the exponential in its Taylor series. The expression
for $g$ (equation (\ref{Yn})) follows as we take the limit
$n\to\infty$ and perform the Fourier transform of (\ref{exampcalc}).

\section{Rate of convergence to stable law}
\label{sconv} 
Since some of the quantities we are trying to estimate depend on 
$n$, it is not enough for us to know that the distributions of 
the quantities $X_n$ converge to a stable law, but it is 
necessary to have good estimates on the speed of convergence (cf. 
Section \ref{limharm})\footnote{Such a result is claimed in 
\cite{kus}, where the authors find an estimate $g_n(x) = g(x) (1 
+ O(1/n))$, where $g_n$ is the density of $X_n$, while $g$ is the 
stable density, and the implicit constants are uniform in $x$. 
Such an estimate would be good enough for our purposes, 
unfortunately it is incorrect, as can be seen by noting that 
$g_n(x) = 0,\quad x< 1-\log n$. A correct estimate is given by 
\cite{phall}, however the method we outline is self-contained, 
simple, and generalizable, so we choose to present it here.}. Our 
results can be summarized as follows:

\thm{theorem} {\label{convthm} In the case $\alpha=1$, we have 
$g_n(x) = g(x) + O(\log^2 n/n).$
}

\begin{proof}
The proof falls naturally in two parts, both of which require
estimates on the characteristic function of $G_n$. First, we show that
we can throw away the tails of characteristic function, and then we
estimate convergence in the remaining region. We introduce the
following notation: let $\Psi_n(k)$ be the characteristic function of
$G_n(k)$, and let $g_n^\delta$ be defined as
\begin{equation*}
g_n^\delta(x) = \frac{1}{2\pi} \int_{-n^{1-\delta}}^{n^{1-\delta}}
\exp(-i x k)\Psi_n(k) d k.
\end{equation*}
In addition, let $R_n(x) = |g_n(x) - g_n^\delta(x)|.$ Then
\thm{lemma} {\label{tailemma}
\begin{equation*}
\sup_{x\in (-\infty, \infty)} R_n(x) = O(ne^{-\frac{\pi}{2}n^{1-\delta}}).
\end{equation*}
}
Lemma \ref{tailemma} will be proved in Section \ref{tail}.

We will also need the following
\thm{lemma}{\label{taillim} Let $\Psi(k)$ be the characteristic
function of the stable density, $g$, and 
\begin{equation*}
R(x)=\frac{1}{2\pi}\left(\int_{-\infty}^{-M} e^{- i k x}
\Psi(k)\,d k + \int_M^\infty e^{- i k x} \Psi(k)\,d k \right). 
\end{equation*}
Then $R(x)= O(\exp(-\pi M/2)).$ }
\begin{proof}
Let us recall that $|\Psi(k)| = \exp(-\frac{\pi}{2} |k|).$ Now we
write
\begin{equation*}
R(x) \leq \frac{1}{2\pi} \left| \int_{-\infty}^{-M} 
e^{-\frac{\pi}{2} |k|}\,d k  +\int_M^\infty e^{-\frac{\pi}{2} |k|}\,
dk \right|= \frac{2}{\pi^2} e^{-\frac{\pi}{2} M}.
\end{equation*}
\end{proof}

To end the proof of Theorem \ref{convthm} we need to estimate how
closely $\Psi(k)$ is approximated by $\Psi_n(k)$, for $k \leq
n^{1-\delta}$, since by the above estimates we know that 
\begin{equation}
\label{denseq}
 g(x) - g_n(x) = 
\frac{1}{\pi}\int_0^{n^{1-\delta}}(\cos(x k) (\Psi(k) - \Psi_n(k))\,d 
k
\end{equation}
(plus lower order terms). It remains only to estimate the difference
$|\Psi(k) - \Psi_n(k)|.$ From equation~(\ref{exampcalc}) we have
\begin{equation*}
 g(x) - g_n(x) \le O(\log^2n/n)\int_0^{n^{1-\delta}}e^{-\frac{\pi}{2}k}k^2\,dk=O(\log^2n/n).
\end{equation*}
\end{proof}

\subsection{\label{tail}An estimate on the tails of the
characteristic function of $Y_n$}
In this section we will supply the proof of Lemma \ref{tailemma};
we use the notation introduced before the statement of the lemma.
We write explicitly
\begin{equation*}
R_n(y)=\Re \int_{n^{1-\delta}}^\infty e^{-ik(y+\log n)}\chi^n(k/n)\,dk.
\end{equation*}
The argument given below can be easily generalized for any stable
law. The only facts about the function $\chi(k)$ that we are going to
use are as follows:
\begin{itemize}
\item[i)]{ $|\chi(k)|\leq 1-c|k|^\beta+o(k)$, $k\ll 1$, for
some $\beta > 0$;}
\item[ii)]{ $|\chi(k)|\le C_1<1$, for all $k>k_0$. This
holds whenever $\chi$ is not the characteristic function of a lattice
distribution, essentially by the Riemann-Lebesgue lemma;}
\item[iii)]{ $\chi \in L^1$.}
\end{itemize}

Given condition (iii) we have

\thm{lemma} {\label{tail3} There exists an $M$ such that
$$\int_M^\infty |\chi(z)|^n  d z < \frac{1}{2^n} 2 \int_M^\infty  |\chi(z)| d
z.$$}
\begin{proof}
By the Riemann-Lebesgue Lemma, $\exists M_2,$ such that $|\chi(z)| \leq
1/2,$ for all $z \geq M_2$. Setting $M=M_2$, the assertion of the
lemma follows immediately.
\end{proof}
Let us introduce the variable $z=k/n$. We have
\begin{equation*}
R_n(y)=n \Re \int_{n^{-\delta}}^\infty e^{-inz(y+\log n)}\chi^n(z)\,dz.
\end{equation*}
Next, we note that
\begin{equation*}
|R_n(y)|\le \tilde{R}_n = n \int_{n^{-\delta}}^\infty 
|\chi^n(z)|\,dz.
\end{equation*}
where the last inequality holds for $\delta < 1$ and $n$ sufficiently
large. Let us write
\begin{equation*}
\tilde{R}_n=R_{n,1}+R_{n,2}+R_{n,3},
\end{equation*}
where
\begin{eqnarray*}
R_{n,1} &=& n\int_{n^{-\delta}}^{z_0}|\chi(z)|^n\,dz,\\
R_{n,2} &=& n\int_{z_0}^{z_1} |\chi(z)|^n\,dz,\\
R_{n,3} &=& n\int_{z_1}^\infty |\chi(z)|^n\,dz.
\end{eqnarray*}
The constant $z_0$ is chosen so that $|\chi(z)|$ is decreasing from
$0$ to $z_0$. Such a $z_0>0 $ can always be found as long as
$\chi(z)$ is continuous at $0$, since $\chi(0) = 1$ and $|\chi(z)| < 1$ for
all $z$ sufficiently close to $0$ (and \textit{all} $z$ if $\chi$ is
not the characteristic function of a lattice distribution). If
$|\chi(z)|$ is monotonically decreasing always, as seems to be the
case with (\ref{concr}), we take $z_0=1$. We choose $z_1$ in such
a way that Lemma \ref{tail3} holds for $M = z_1$.

First we consider $R_{n,1}$. By the choice of $z_0$, the function
$|\chi(z)|$ monotonically decreases on $[n^{-\delta},z_0]$. Therefore,
\begin{equation}
\label{R1}
|R_{n,1}|\le nz_0|\chi(n^{-\delta})|^n.
\end{equation}
For small values of $k$, we have (see property (i)):
\begin{equation*}
|\chi(k)|=1-\pi k/2+O((k\log k)^2).
\end{equation*}
Therefore,
\begin{equation*}
|\chi(n^{-\delta})|^n\asymp \exp(-\frac{\pi}{2} n^{1-\delta}),
\end{equation*}
that is, it decays exponentially for any $0<\delta<1$. Thus, from
(\ref{R1}) we see that
\begin{equation}
\label{R10} |R_{n,1}|\le nz_0 \exp(-\frac{\pi}{2} n^{1-\delta}).
\end{equation}
Next, we estimate $R_{n,2}$. Because of property (ii), we have
\begin{equation}
\label{R2} |R_{n,2}|\le n^{1+\delta}C_1^n,
\end{equation}
and since $C_1 < 1$, it also decays exponentially with $n$.  Finally,
we estimate $R_{n,3}$ by Lemma \ref{tail3}.  Putting everything
together we conclude that $R_n$ decays exponentially with $n$.

\end{document}